\documentclass[12pt]{amsart}
\usepackage[colorlinks=true,citecolor=blue,linkcolor=blue]{hyperref}
\usepackage{amssymb,amscd,amsmath,graphicx,mathtools}
\usepackage{enumerate}
\usepackage{fullpage}
\usepackage{comment}
\usepackage{xcolor}
\usepackage[all]{xy}
\usepackage{cleveref}
\usepackage{enumerate}
\usepackage{fullpage}
\usepackage{xcolor}
\usepackage{amsthm}
\usepackage{float}
\usepackage{pgfplots}
\usepackage{listings}
\usepackage{bbm}
\usepackage{comment}
\usepackage{graphicx}

\numberwithin{equation}{section}
\newtheorem{theorem}{Theorem}[section]
\newtheorem{lemma}[theorem]{Lemma}

\newtheorem{corollary}[theorem]{Corollary}

\theoremstyle{definition}
\newtheorem{definition}[theorem]{Definition}

\newtheorem{question}[theorem]{Question}

\theoremstyle{definition}
\newtheorem{remark}[theorem]{Remark}
\newtheorem{example}[theorem]{Example}

\DeclareMathOperator{\supp}{supp}

\DeclareMathOperator{\NP}{NP}

\DeclareMathOperator{\conv}{convexhull}

\newcommand{\CC}{{\mathbb C}}

\newcommand{\ZZ}{{\mathbb Z}}
\newcommand{\NN}{{\mathbb N}}

\newcommand{\RR}{{\mathbb R}}

\def\I{{\mathcal I}}

\def\C{{\mathcal C}}

\def\R{{\mathcal R}}

\def\R{{\mathcal R}}

\def\mm{{\mathfrak m}}

\def\pp{{\mathfrak p}}

\def\a{{\bf a}}
\def\b{{\bf b}}

\def\x{{\bf x}}

\newcommand{\kk}{\Bbbk}

\def\1{{\bf 1}}
\def\0{{\bf 0}}

\begin{document}
\title[Multiplicity = Volume formula in regular local rings]{Multiplicity = Volume formula and Newton non-degenerate ideals in regular local rings}

\author{T\`ai Huy H\`a}
\address{Tulane University, Department of Mathematics, 6823 St. Charles Ave., New Orleans, LA 70118, USA}
\email{tha@tulane.edu}

\author{Th\'ai Th\`anh Nguy$\tilde{\text{\^E}}$n}
\address{University of Dayton, Department of Mathematics,
	300 College Park, Dayton, Ohio, USA \\
	and University of Education, Hue University, 34 Le Loi, Hue, Vietnam}
\email{tnguyen5@udayton.edu}

\author{Vinh Anh Ph{\d{A}}m}
\address{Tulane University, Department of Mathematics, 6823 St. Charles Ave., New Orleans, LA 70118, USA}
\email{vpham1@tulane.edu}






\begin{abstract}
We develop the notions of Newton non-degenerate (NND) ideals and Newton polyhedra for regular local rings. These concepts were first defined in the context of complex analysis. We show that the characterization of NND ideals via their integral closures known in the analytical setting extends to regular local rings. We use the limiting body $\C(\I)$ associated to a graded family $\I$ of ideals to provide a new understanding of the celebrated ``Multiplicity $=$ Volume'' formula. Particularly, we prove that, for a Noetherian graded family $\I$ of $\mm$-primary ideals in a regular local ring $(R,\mm)$ of dimension $d$, the equality
$$e(\I) = d!\text{co-vol}_d(\C(\I))$$
holds if and only if $\I$ contains certain subfamily of NND ideals.
\end{abstract}

\maketitle

\section{Introduction} \label{sec.intro}

Let $(R,\mm)$ be a regular local ring of dimension $d$ containing a field, and let $\I = \{I_n\}_{n \in \NN}$ be a graded family of $\mm$-primary ideals in $R$. A well-celebrated result in the theory of Newton-Okounkov bodies, that was developed by Lazarsfeld--Musta\c{t}\u{a} \cite{LM09}, Kaveh--Khovanskii \cite{KK12, KK14} and Cutkosky \cite{Cut13, Cut14}, and extended in various directions by Cid-Ruiz--Mohammadi--Monin \cite{CMM24}, Cid-Ruiz--Monta\~no \cite{CJ22a, CJ22b}, is the ``Multiplicity $=$ Volume'' formula. This formula generalizes Teissier's classical formula for monomial ideals given in \cite{Tei88}, and states that the \emph{multiplicity} $e(\I) = \displaystyle \limsup_{n \rightarrow \infty} \frac{\ell_R(R/I_n)}{n^d/d!}$ exists as an actual limit and can be computed from the integral volume (in $\RR^d$) of the \emph{Newton-Okounkov body} associated to $\I$.
In general, the construction of Newton-Okounkov bodies relies on appropriate \emph{good} valuations, which are often challenging to understand. The motivating question for our work is whether and when other convex bodies associated to $\I$ could replace the Newton-Okounkov body to offer a more accessible approach to the ``Multiplicity $=$ Volume'' formula. Our results show that such substitutions are indeed feasible when the family $\I$ contains a subfamily of \emph{Newton non-degenerate} ideals.

Newton non-degenerate ideals are interesting and have their own important motivations. The concepts of Newton non-degenerate ideals and their associated \emph{Newton polyhedra} were initially introduced in the context of complex analysis by Saia \cite{Saia96} to study geometric invariant in the ring of germs of holomorphic functions at the origin of $\CC^d$. It was later extended by Bivi\`a-Ausina \emph{et. al.} \cite{BA04, BAFS02} to formal power series rings over the complex numbers. Newton polyhedra of ideals generalize the more familiar notion for \emph{monomial} ideals in polynomial rings. Newton non-degenerated ideals and Newton polyhedra also appeared in several research areas of commutative algebra and algebraic combinatorics, including the study of core of ideals \cite{PUV07}, torus-closure of algebraic schemes \cite{Chen17} and in exploring symmetric polynomials \cite{Stan84, Stem91}, Schur polynomial \cite{Rado52}, symmetric Macdonald and Schubert polynomials \cite{MTY19}.

An essential fact in understanding Newton non-degenerate ideals is that, for the ring of germs of holomorphic functions at the origin or formal power series rings, these ideals are characterized by the property that their integral closures are monomial ideals (see \cite{BAFS02, Saia96}). However, this characterization does not carry over to, for example, polynomial rings (see \Cref{rmk.NonLocal}).
To address our motivating question about the ``Multiplicity $=$ Volume'' formula, we develop the notions of Newton non-degenerate ideals and Newton polyhedra for regular local rings. We also demonstrate that the characterization of Newton non-degenerate ideals via their integral closures known in the analytical setting extends to regular local rings.

More specifically, let $\x = x_1, \dots, x_d$ be a regular system of parameters in $R$. A \emph{monomial} in $\x$ is an element in $R$ of the form $\x^\alpha = x_1^{\alpha_1} \dots x_d^{\alpha_d}$, where $\alpha = (\alpha_1, \dots, \alpha_d) \in \ZZ_{\ge 0}^d$.
For a nonzero element $f \in R$, write $f = \sum_{i=1}^t r_i \cdot \x^{\a_i}$, in which $r_i \in R \setminus \mm$ for all $i$, and set
$\supp (f) = \{\a_1, \dots, \a_t\}.$
Roughly speaking, the \emph{Newton polyhedron} of an ideal $I \subseteq R$ (with respect to $\x$) is
$$\Gamma_\x(I) = \conv \langle \left\{u ~\big|~ u \in \supp (f), f \in I\right\}\rangle \subseteq \RR_{\ge 0}^d.$$
We show that Newton polyhedra of ideals in regular local rings exhibit many important properties as those known in the analytic case. For example, we prove that the Newton polyhedron $\Gamma_\x(I)$ of an ideal $I \subseteq R$ can be obtained from those of a generating set for $I$ and does not depend on the particular choices of the generators; see \Cref{thm.NPgenerators}.





For each face $\Delta \subseteq \Gamma_+(I)$, the set of points on rays through $\Delta$ emanating from the origin $\0 = (0, \dots, 0) \in \RR^d$ forms a cone $C(\Delta)$. The intersection $C(\Delta) \cap \ZZ_{\ge 0}^d$ is a subsemigroup of $\ZZ^d$, yielding a local ring
$$R_\Delta = \{g \in R ~\big|~ \supp (g) \subseteq C(\Delta) \cap \ZZ_{\ge 0}^d\},$$
a subring of $R$ with a unique maximal ideal
$$\mm_\Delta = \{f \in R_\Delta ~\big|~ \0 = (0, \dots, 0) \not\in \supp (f)\}.$$
Given $f = \sum_{i=1}^t r_i \x^{\a_i} \in R$, with $r_i \in R \setminus \mm$, set $f_\Delta = \sum_{\a_i \in C(\Delta)} r_i\x^{\a_i}.$ An ideal $I \subseteq R$ is call a \emph{Newton non-degenerate} ideal if there exists a system of generators $\{g_1, \dots, g_s\}$ of $I$ such that, for each compact fact $\Delta \subseteq \Gamma_+(I)$, the ideal $I_\Delta = \left( g_{1\Delta}, \dots, g_{s\Delta}\right)$ is $\mm_\Delta$-primary in $R_\Delta$.

Our results extend the characterization of Newton non-degenerate ideals via their integral closures to regular local rings. We prove the following theorem.

\medskip

\noindent\textbf{Theorem \ref{thm.NNDequiv}.} Let $(R,\mm)$ be a regular local ring of dimension $d$ with a regular system of parameters $\x = x_1, \dots, x_d$. Let $I \subseteq R$ be an ideal. Then, $I$ is Newton non-degenerate if and only if its integral closure $\overline{I}$ is a monomial ideal in $\x$.

\medskip

To establish \Cref{thm.NNDequiv}, we consider the faithfully flat extension $R \rightarrow \widehat{R}$ and show that $I_\Delta$ is $\mm_\Delta$-primary if and only if $(I\widehat{R})_\Delta$ is $\widehat{\mm}_\Delta$-primary; see \Cref{cor_m_primary_between_R_and_hatR}. Since $\widehat{R}$ is isomorphic to a power series ring over the field $k \simeq R/\mm$, the known result in this case then kicks in.

The concepts of Newton non-degenerate ideals and their Newton polyhedra facilitate our study of the ``Multiplicity $=$ Volume'' formula. Particularly, for a graded family $\I = \{I_n\}_{n \in \NN}$, set
$$\C(\I) = \bigcup_{n \in \NN} \frac{1}{n}\Gamma_\x(I_n) \subseteq \RR_{\ge 0}^d.$$
The set $\C(\I)$ is usually referred to as the \emph{limiting body} of $\I$. This notion has appeared in other contexts (eg, cf. \cite{May14a, Wol08}). It was also introduced and studied for graded families of monomial ideals in polynomial rings (see \cite{CDF+, HN24}).
Our results demonstrate that the limiting body can replace the Newton-Okounkov body in the ``Multiplicity $=$ Volume'' formula when the family $\I$ contains a subfamily of Newton non-degenerate ideals. We establish the following theorem.

\medskip

\noindent\textbf{Theorem \ref{thm.Mult}.} Let $(R,\mm)$ be a regular local ring of dimension $d$ and let $\x = x_1, \dots, x_d$ be a regular system of parameters in $R$. Let $\I = \{I_n\}_{n \in \NN}$ be a Noetherian graded family of $\mm$-primary ideals in $R$. Let $c$ be an integer such that $\overline{I_c^k} = \overline{I_{kc}}$, for all $k \in \NN$, as in \Cref{thm.Noetherian}. Then, the following are equivalent:
\begin{enumerate}
\item $e(\I) = d! \text{co-vol}_d (\C(\I))$, and
\item $I_{kc}$ is an NND ideal for every $k \in \NN$.
\end{enumerate}

\medskip

\Cref{thm.Mult} is achieved by considering the ideal $I_0$ associated to each ideal $I \subseteq R$, that is generated by monomials $\{\x^\a ~\big|~ \a \in \Gamma_\x(I) \cap \ZZ_{\ge 0}^d\}$, and understanding the relationship between $I_0$ and $\overline{I}$, as well as that between $e(I_0)$ and $e(I)$, via $\Gamma_\x(I)$; see \Cref{thm.NNDregular} and \Cref{cor.IntegralClosureofNNDidealsPoly}.

Our paper is outlined as follows. In \Cref{sec.Monomials}, we develop the theory of monomials, monomial ideals and Newton polyhedra of arbitrary ideals for regular domains; these are domains whose localization at any prime is regular. We show that these concepts share many similarities with the more familiar ones in polynomial rings or rings of formal power series. In \Cref{sec.NND}, we generalize the notion of Newton non-degenerate ideals to regular local rings and illustrate that the characterization of these ideals via their integral closures holds. In \Cref{sec.Multiplicity}, we consider graded families of ideals and the ``Multiplicity $=$ Volume'' formula. We also obtain a classification for a graded family of Newton non-degenerate ideals to be Noetherian.


\medskip

\noindent\textbf{Acknowledgment.} The first author is partially supported by a Simons Foundation grant. This is part of the third author's PhD thesis. 


\section{Newton polyhedra associated to ideals in regular domains} \label{sec.Monomials}

In this section, we develop the notion of Newton polyhedra associated to ideals in \emph{regular} domains; these are integral domains whose localization at every prime ideal is regular. This is possible partly thanks to the following notion of generalized regular systems of parameters.

\begin{definition}[{\cite[Definition 1.1]{HS08}}] \label{def.genSOP}
Let $R$ be a regular domain. Elements $x_1, \dots, x_p$ in $R$ are called a \emph{generalized regular system of parameters} (g.r.s.o.p) if $x_1, \dots, x_p$ is a permutable regular sequence in $R$ such that, for any collection $i_1, \dots, i_t \subseteq \{1, \dots, p\}$, $R/(x_{i_1}, \dots, x_{i_t})$ is a regular domain.
\end{definition}

It is easy to see that any part of a g.r.s.o.p is a g.r.s.o.p. It was also remarked in \cite{HS08} that when $R$ is a regular local ring, an arbitrary \emph{regular system of parameters} (r.s.o.p) (or a part thereof) is a generalized regular system of parameters, and when $R$ is a polynomial ring over a field, the variables form a generalized regular system of parameters.

\medskip

\noindent\textbf{Set-up.} Throughout this section, $R$ is a regular domain and $\x = x_1, \dots, x_p$ denotes a fixed g.r.s.o.p in $R$.

\medskip

The concepts of g.r.s.o.p enables us to define ``monomials'', ``monomial ideals'' and ``Newton polyhedron'' in a manner similar to the more familiar case of a polynomial ring. However,  special care is required at times, because, for instance, an ideal $I$ being monomial is not necessarily equivalent to the fact that $f \in I$ if and only if all monomials appearing in $f$ are in $I$. In the following, for $\a= (a_1,\ldots ,a_p), \b=(b_1,\ldots ,b_p)\in \RR^p$, by $\a \ge \b$, we mean $a_i\ge b_i$ for all $i=1,\ldots,p$.

\begin{definition}[See {\cite[Definitions 2.1 and 2.2]{HS08}}] \label{def.monomial}
Let $R$ be a regular domain with a g.r.s.o.p $\x = x_1, \dots, x_p$.
\begin{enumerate}
    \item A \emph{monomial} in $\x$ is an element in $R$ of the form $\x^\alpha = x_1^{\alpha_1} \cdots x_p^{\alpha_p}$, where $\alpha = (\alpha_1, \dots, \alpha_p) \in \ZZ_{\ge 0}^p$.
    \item A \emph{monomial ideal} in $R$ with respect to $\x$ is an ideal generated by monomials in $\x$.
    \item If $I = (\x^{\a_1}, \dots, \x^{\a_s})$ is a monomial ideal in $R$ with respect to $\x$, then the \emph{Newton polyhedron} of $I$ is defined as
        $$\NP(I) = \left\{ \a \in \RR_{\ge 0}^p ~\Big|~ \a \ge \sum_{i=1}^s c_i \a_i \text{ for some } c_i \in \RR_{\ge 0} \text{ such that } \sum_{i=1}^s c_i = 1\right\} \subseteq \RR_{\ge 0}^p.$$
\end{enumerate}
\end{definition}

\begin{example} \label{ex.1}
Let $R=\mathbb{R}[x,y]$. Consider the sequence $\x = x^2+1,y^2+2$. It is not hard to check that this sequence is a g.r.s.o.p. Consider an ideal $I = (x^2y^2+2x^2+y^2+2, y^4+4y^2+4)$ in $R$. Since
$I=((x^2+1)(y^2+2), (y^2+2)^2)$, then $\NP(I)$ will be the convex hull $\{(1,1), (0,2)\}+\RR^2_{\geq 0}$.
\begin{figure}[h]
	\caption{$\NP(I)$ for $I = ((x^2+1)(y^2+2), (y^2+2)^2) \subseteq \RR[x,y]$.}
	\centering
	\includegraphics[width=0.4\textwidth]{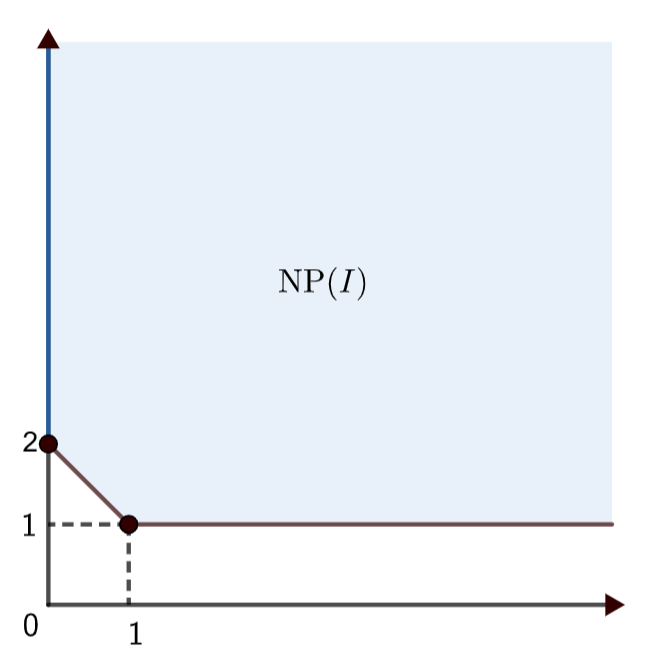}
	\label{NP_monomial_ideal}
\end{figure}
\end{example}

It is easy to see that $\NP(I)$, for a monomial ideal $I$, is a closed convex set. The following theorem, proved in \cite{HS08}, allows us to express elements in $R$ in terms of its monomials. Note that \cite[Theorem 1.3]{HS08} was stated for nonzero element $f \in R$ and the statement is trivial when $f = 0$.

\begin{theorem}[{\cite[Theorem 1.3]{HS08}}] \label{thm.HS08}
Let $\pp = (\x) = (x_1,\dots,x_p)$ and let $f \in R$. Then, there exist monomials $m_1,\dots,m_t$ in $\x$ and elements $h, r_1,\dots, r_t\in R\setminus \pp$ such that
		$$h\cdot f=\sum\limits_{i=1}^tr_i\cdot m_i.$$
\end{theorem}

We call the expression $h\cdot f = \sum_{i=1}^t r_i \cdot m_i$ as in \Cref{thm.HS08} a \emph{monomial representation} of $f$ with respect to $\x$. The following observations show that from a monomial representation of an element in $R$, we can always reduce it to an \emph{irredundant} one, and furthermore, this irredundant monomial representation is "unique".

\begin{corollary} \label{cor.HS08irred}
    Let $f \in R$ and let
$$h \cdot f = \sum_{i=1}^t r_i \cdot m_i,$$
where $m_1, \dots, m_t$ are monomials in $\x$, $h, r_1, \dots, r_t \in R \setminus \pp$, be a monomial representation of $f$ as in \Cref{thm.HS08}. Then, there exists a subset $\{m_{i_1},\ldots ,m_{i_k} \} \subseteq \{m_1,\ldots,m_t \}$ and elements $r_1',\ldots ,r_k' \in R \setminus \pp$ such that $m_{i_j} \nmid m_{i_\ell}$ for any $j\ne \ell, j, \ell =1,\ldots, k$ such that $$h \cdot f = \sum_{j=1}^k r_j' \cdot m_{i_j}.$$
We call this an \emph{irredundant} monomial representation of $f$.
\end{corollary}
\begin{proof}
    Suppose there exists $i,j$ such that $m_i \ | \ m_j$, then we can write $m_j=m_i' \cdot m_i$ where $m_i'$ is a monomial. It follows that we can write
\begin{align*}
    hf&= \sum_{l\neq i,j}^t r_l m_l +r_im_i+r_jm_j= \sum_{l\neq i,j}^t r_l m_l +(r_i+r_jm_i')m_i.
\end{align*}
The element $r_i+r_jm_i'$ is not in $\mathfrak{p}$ since otherwise $r_i\in \mathfrak{p}$ (a contradiction). It means that we have another monomial representation of $f$ without $m_j$. Repeat this process to remove from the monomial presentation any $m_j$ such that there is an $m_i$ dividing it, since the set of monomials in a representation is finite, this process must terminate in finite steps. At the end, we obtain a monomial representation of $f$ such that $m_i\nmid m_j$ for all $i,j$.
\end{proof}

We call the above monomial representation an \textit{irredundant} one.

\begin{corollary} \label{cor.HS08unique}
Let $f \in R$ and let
$$h \cdot f = \sum_{i=1}^t r_i \cdot m_i,$$
where $m_1, \dots, m_t$ are monomials in $\x$, $h, r_1, \dots, r_t \in R \setminus \pp$, be a monomial representation of $f$ as in \Cref{thm.HS08}. Suppose, in addition, that $m_i \nmid m_k$ for any $i \not= k$. Then, the set $\{m_1, \dots, m_t\}$ is unique. We also call the above a \emph{minimal} monomial representation of $f$.
\end{corollary}

\begin{proof}
By localizing $R$ at $\pp$, we may assume that $R$ is a regular local ring with maximal ideal $\pp$. Suppose that we have two monomial representations of $f$, namely,
	\begin{equation*}
		h\cdot f=\sum\limits_{i=1}^tr_i\cdot m_i \text{\space\space and \space\space} h'\cdot f=\sum\limits_{j=1}^rt_j\cdot n_j,
	\end{equation*}
	where $m_i\nmid m_k$ for all $i,k\in \{1,\dots,t\}$ with $i\neq k$, and $n_j\nmid n_l$ for all $j,l\in \{1,\dots,r\}$ with $j\neq l$. Then,
	\begin{equation*}
		hh'\cdot f=\sum\limits_{i=1}^tr'_i\cdot m_i=\sum\limits_{j=1}^rt'_j\cdot n_j,
	\end{equation*}
	where $r_i'=h'r_i$ for $i=1,\dots,t$, and $t'_j=ht_j$ for $j=1\dots,r$.

Since $r_i'$ and $t'_j$ are not in the maximal ideal $\pp$ for $i=1,\dots,t$ and for $j=1\dots,r$, then they must be units in $R$. It follows that for each $i \in \{1,\dots,t\}$, we can write
\[
m_i = \sum\limits_{j=1}^r a_j\cdot n_j + \sum\limits_{j\ne i}^t b_j\cdot m_j,
\]
for some units $a_j,b_j$ in $R$. As, $m_j\nmid m_i$ for all $j \ne i$, by \cite[Corollary 3]{KS03}, there exist $j\in \{1,\dots,r\}$ such that $n_j\mid m_i$. Likewise, there exits $k\in \{1,\dots,t\}$ with such that $m_k\mid n_j$. This implies that $m_k\mid m_i$. By the irredundance hypothesis, this is possible only if $i=k$ and $m_i=n_j$. In particular, it follows that $\{m_1,\dots,m_t\}\subseteq \{n_1,\dots, n_r\}$. A similar argument shows that $\{n_1,\dots, n_r\}\subseteq \{m_1,\dots,m_t\}$. The assertion is proved.
\end{proof}

In light of \Cref{cor.HS08unique}, from now on, when we refer to a monomial representation of an element in $R$, we mean its unique irredundant monomial representation. \Cref{cor.HS08unique} allows us to define the support of an element in $R$ with respect to a g.r.s.o.p.

\begin{definition} \label{def.supp}
Let $R$ be a regular domain and let $\x = x_1, \dots, x_p$ be a g.r.s.o.p in $R$.
\begin{enumerate}
\item Let $f \in R$ and suppose that the unique monomial representation of $f$ with respect to $\x$ is
$$h \cdot f = \sum_{i=1}^t r_i \cdot m_i,$$
where $m_1 = \x^{\a_1}, \dots, m_t = \x^{\a_t}$ are monomials in $f$ and $h, r_1,\ldots,r_t \not\in (\x)=(x_1,\ldots,x_p)$. The \emph{support} of $f$ is defined as
$$\supp(f) = \{\a_1, \dots, \a_t\} \subseteq \RR_{\ge 0}^p.$$
\item Let $I \subseteq R$ be an ideal. The \emph{support} of $I$ is defined to be
    $$\supp(I) = \bigcup_{f \in I} \supp(f) \subseteq \RR_{\ge 0}^p.$$
\end{enumerate}
\end{definition}

We are now ready to extend the notion of Newton polyhedra from monomial ideals to any ideal in $R$. This definition generalizes that initially given in complex analysis (cf. \cite{Saia96}).

\begin{definition} \label{def.Newton}
Let $R$ be a regular domain with a g.r.s.o.p $\x = x_1, \dots, x_p$. Let $I \subseteq R$ be an ideal. The \emph{Newton polyhedron} of $I$ (w.r.t. $\x$) is
$$\Gamma_{\x}(I) = \conv\langle \left\{\a \in \ZZ^p_{\ge 0} ~\big|~ \a \in \supp(I)\right\}\rangle \subseteq \RR_{\ge 0}^p.$$
\end{definition}


\begin{example} \label{ex.2}
Let $R=\mathbb{R}[x,y]$. Consider the sequence $\x = x^2+1,y^2+2$. From Example \ref{ex.1}, this sequence is a g.r.s.o.p. Consider $f = (x^2+1)^2(y^2+2)+(x^2+1)(y^2+2)^3$ and the ideal $I = (f)$. Then, the support of $I$ is $\supp(I)=\{(2,1), (1,3)\}$. Therefore, the Newton polyhedron $\Gamma_\x(I)$ will be the convex hull $\{(2,1), (1,3)\}+\RR^2_{\geq 0}$, as depicted in \Cref{Gamma_x_ideal}.
\begin{figure}[h]
	\caption{$\Gamma_\x(I)$ for $I = ((x^2+1)^2(y^2+2)+(x^2+1)(y^2+2)^3) \subseteq \RR[x,y]$.}
	\centering
	\includegraphics[width=0.4\textwidth]{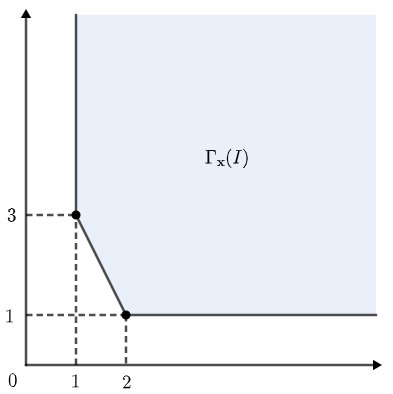}
	\label{Gamma_x_ideal}
\end{figure}
\end{example}

\begin{remark} \label{rmk.NPPrincipal} If $I = (f)$ is a principal ideal, then we denote by $\Gamma_\x(f)$ the Newton polyhedron of $I$. It follows from the definition that
$$\Gamma_\x(f) = \conv\left\langle \{\supp(f) + \supp(h) ~\big|~ h \in R\}\right\rangle = \conv\langle \supp(f)\rangle + \RR_{\ge 0}^p.$$
\end{remark}

The next lemma gives basic properties of supports.

\begin{lemma}\label{lem.properties.of.supp}
	Let $f,g$ be elements in $(\x) = (x_1, \dots, x_p)$. Then,
	\begin{enumerate}
		\item $\supp(fg)\subseteq \{\a + \b ~\big|~ \a \in \supp(f), \b \in \supp(g)\} = \supp(f)+\supp(g)$; and
		\item $\supp(f+g)\subseteq \supp(f)\cup \supp(g)$.
	\end{enumerate}
\end{lemma}

\begin{proof} Suppose that the monomial representations of $f$ and $g$ with respect to $\x$ are
$$h \cdot f = \sum_{i=1}^t r_i m_i \text{ and } h' \cdot g = \sum_{j=1}^r r_j' m_j',$$
where $m_1, \dots, m_t, m_1', \dots, m_r'$ are monomials in $\x$, $h, h', r_i, r_j' \in R \setminus \pp$ with $\pp = (\x)$.

(1) Since $h, h' \not\in \pp$, $hh' \not\in \pp$. Thus,
$$(hh') \cdot (fg) = \left(\sum_{i=1}^t r_i m_i\right) \cdot \left(\sum_{j=1}^r r_j' m_j'\right)$$
is a monomial representation of $fg$. By Corollaries \ref{cor.HS08irred}, \ref{cor.HS08unique}, monomials appearing in a minimal monomial representation of $fg$ are of the form $m_im_j'$ for some $i,j$. This implies that
$$\supp(fg) \subseteq \{\a + \b ~\big|~ \a \in \supp(f), \b \in \supp(g)\} = \supp(f) + \supp(g).$$

(2) We have
$$(hh')\cdot (f+g) = \sum_{i=1}^t (h'r_i)m_i + \sum_{j=1}^r (hr_j')m_j'.$$
Since $h,h', r_i, r_j' \not\in \pp$, $h'r_i \not\in \pp$ and $hr_j' \not\in \pp$ for any $i,j$. Again, by Corollaries \ref{cor.HS08irred}, \ref{cor.HS08unique}, monomials in a minimal representation of $f+g$ are among those $m_i,m_j'$'s. It follows that
$$\supp(f+g) \subseteq \supp(f) \cup \supp(g). \qedhere$$
\end{proof}

Our first result shows that, like the case of rings of germs of holomorphic functions or formal power series rings, $\Gamma_{\x}(I)$ can be constructed from a set of generators of $I$ and does not depend on the choices of the generators.

\begin{theorem}\label{thm.NPgenerators}
	Let $I \subseteq R$ be an ideal and let $g_1,\dots,g_s$ be a generating set of $I$. Then,
	\begin{equation*}
		\Gamma_\x(I)=\conv\left\langle\Gamma_\x(g_1)\cup\cdots\cup\Gamma_\x(g_s)\right\rangle.
	\end{equation*}
\end{theorem}

\begin{proof}
	For each $i=1,\dots,s$, $g_i \in I$, so $\Gamma_\x(g_i)\subseteq \Gamma_\x(I)$. Since $\Gamma_\x(I)$ is convex, it follows that
$$\conv\langle \Gamma_\x(g_1)\cup\cdots\cup\Gamma_\x(g_s) \rangle\subseteq \Gamma_\x(I).$$
	
	Conversely, let $f$ be an element of $I$. It follows that $f=\sum_{i=1}^s h_ig_i$ for some $h_i\in R$. By Lemma \ref{lem.properties.of.supp}, we have
$$\supp(f) \subseteq \bigcup_{i=1}^s \supp(h_ig_i) \subseteq \bigcup_{i=1}^s [\supp(h_i) + \supp(g_i)] \subseteq \bigcup_{i=1}^s [\supp(g_i) + \RR_{\ge 0}^p].$$
As observed in Remark \ref{rmk.NPPrincipal}, $\conv\langle\supp(g_i)\rangle + \RR_{\ge 0}^p = \Gamma_\x(g_i)$. By letting $f$ run through $I$, we get
$$\Gamma_\x(I)\subseteq \conv\langle\Gamma_\x(g_1)\cup\cdots\cup\Gamma_\x(g_s)\rangle,$$
and the statement is proved.
\end{proof}

\begin{example} \label{ex.3}
	 Consider a regular local ring $(R,\mathfrak{m})$ with $\x=x,y,z$ as a regular system of parameters. By definition, the sequence is also a g.r.s.o.p. Let $I=(x^4+x^2y^5z^4,y^3z-x^3yz^2+xyz^3)$. Then by Theorem \ref{thm.NPgenerators}, the Newton polyhedron of $I$ will be the convex hull of the union of $\Gamma_{\x}(x^4+x^2y^5z^4)$ and $\Gamma_{\x}(y^3z-x^3yz^2+xyz^3)$, as depicted in \Cref{Gamma_by_generators}.
	 \begin{figure}[h]
	 	\caption{$\Gamma_\x(I)$ for $I = (x^4+x^2y^5z^4,y^3z-x^3yz^2+xyz^3) \subseteq R$.}
	 	\centering
	 	\includegraphics[width=0.4\textwidth]{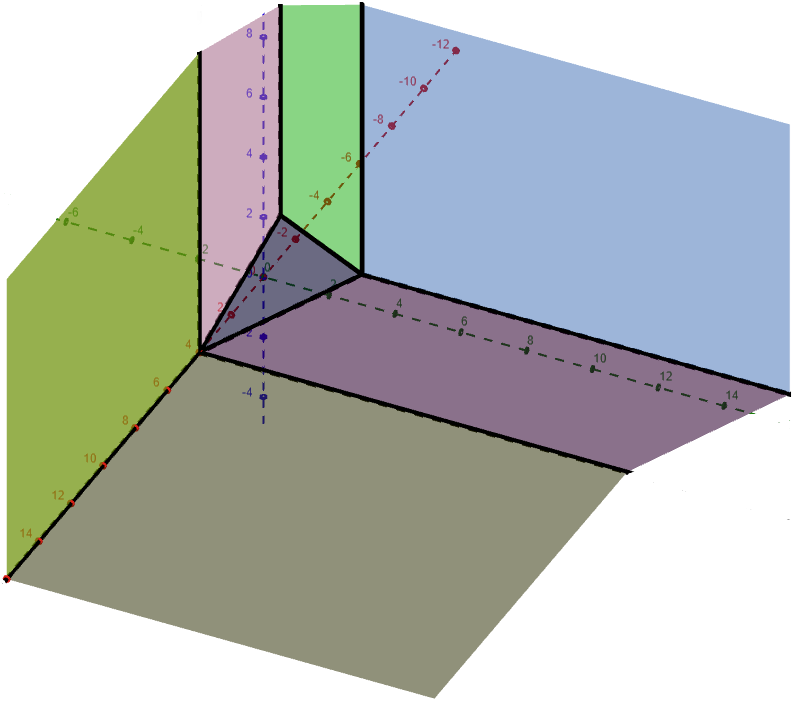}
	 	\label{Gamma_by_generators}
	 \end{figure}
\end{example}

\begin{corollary} \label{cor.Gamma=NP}
If $I \subseteq R$ is a monomial ideal in $\x$, then
$$\Gamma_\x(I) = \NP(I).$$
\end{corollary}

\begin{proof} Suppose that $I = (\x^{\a_1}, \dots, \x^{\a_s})$. By \Cref{thm.NPgenerators} and \Cref{def.monomial}, we have
$$\Gamma_\x(I) = \conv \langle \Gamma_\x(\x^{\a_1}), \dots, \Gamma_\x(\x^{\a_s})\rangle = \conv \left\langle \bigcup_{i=1}^s [\a_i + \RR_{\ge 0}^p]\right\rangle = \NP(I). \qedhere$$
\end{proof}

Recall that the \emph{Minkowski sum} of subsets $A, B \subseteq \RR^p$ is given by
$A+B=\{a+b \mid a\in A, b\in B\}.$ It is also a standard fact that if $A$ and $B$ are convex sets then so is their Minkowski sum $A+B$. The following lemmas are basic properties of Minkowski sums. We include the proofs due to the lack of appropriate references.

\begin{lemma}\label{lem1}
	Let $P, Q \subseteq \RR^p$ be two polyhedra. If $(a,c)$ and $(b,d)$ are different points in $P \times Q$ such that $a+c=b+d$, then $a+c$ is not a vertex of $P+Q$.
\end{lemma}

\begin{proof}
	If $a+c = b+d$, then
	\begin{equation*}
		a+c=\frac{1}{2}(a+d)+\frac{1}{2}(b+c).
	\end{equation*}
	Since $a+d \not= a+c$ and $b+c \not= a+c$ also belong to $P+Q$, $a+c$ is not a vertex of $P+Q$.
\end{proof}

\begin{lemma}\label{lem.Gammaprod=sumGamma}
	Let $g, h \in R$. Then
	\begin{equation*}
		\Gamma_\x(gh)=\Gamma_\x(g)+\Gamma_\x(h)
	\end{equation*}
\end{lemma}

\begin{proof} By \Cref{lem.properties.of.supp}, $\supp(gh)\subseteq \supp(g) +\supp (h).$ Thus, 	$\Gamma_\x(gh)\subseteq\Gamma_\x(g)+\Gamma_\x(h)$.

To establish the reverse inclusion, it is enough to show that every vertex of $\Gamma_\x(g)+\Gamma_\x(h)$ is in $\Gamma_\x(gh)$. Observe that each vertex of $\Gamma_\x(g)+\Gamma_\x(h)$ is of the form $u+v$, where $u$ is a vertex of $\Gamma_\x(g)$ and $v$ is a vertex of $\Gamma_\x(h)$. Particularly, $u\in \supp(g)$ and $v\in \supp(h)$.

If $u+v\notin \supp(gh)$, then the monomial $\x^u \cdot \x^v$ does not appear in $gh$. Therefore, there exist $a\in \supp(g)$ and $b\in \supp(h)$, with $(u,v)$ and $(a,b)$ being different points in $\Gamma_\x(g) \times \Gamma_\x(h)$, such that $u+v=a+b$. This, however, is a contradiction to the conclusion of Lemma \ref{lem1}. Hence, $u+v \in \supp(gh)$ and the statement is proved.
\end{proof}



The next result extends \Cref{lem.Gammaprod=sumGamma} to a more general setting.

\begin{theorem} \label{thm.product}
	Let $I,J\subseteq R$ be ideals. Then,
	\begin{equation*}
		\Gamma_\x(IJ) = \Gamma_\x(I) + \Gamma_\x(J).
	\end{equation*}
\end{theorem}

\begin{proof}
	Suppose that $\{g_1,\dots,g_u\}$ and $\{h_1,\dots,h_v\}$ are generating sets of $I$ and $J$, respectively. Then, $IJ$ is generated by the set $\{g_ih_j \mid 1\leq i\leq u, 1\leq j\leq v\}$. By Theorem \ref{thm.NPgenerators}, we have
	\begin{equation*}
		\Gamma_{\x}(IJ)=\conv \left\langle \bigcup\limits_{1\leq i\leq u,1\leq j\leq v}\Gamma_\x(g_ih_j)\right\rangle.
	\end{equation*}
	For each $i$ and $j$, $\Gamma_\x(g_ih_j)=\Gamma_\x(g_i)+\Gamma_\x(h_j)\subseteq \Gamma_\x(I)+\Gamma_\x(J)$, by Lemma \ref{lem.Gammaprod=sumGamma}. It follows that
	\begin{equation*}
		\bigcup\limits_{1\leq i\leq u,1\leq j\leq v}\Gamma_\x(g_ih_j)\subseteq \Gamma_\x(I)+\Gamma_\x(J),
	\end{equation*}
	and, hence, $\Gamma_{\x}(IJ)\subseteq \Gamma_\x(I)+\Gamma_\x(J)$.
	
	Conversely, we have
	\begin{align*}
		\Gamma_\x(I)+\Gamma_\x(J)&= \conv \left\langle \bigcup\limits_{i=1}^{u}\Gamma_\x(g_i)\right\rangle + \conv \left\langle \bigcup\limits_{j=1}^{v}\Gamma_\x(h_j)\right\rangle\\
		&=\conv  \left\langle \left(\bigcup\limits_{i=1}^{u}\Gamma_\x(g_i)\right)+\left(\bigcup\limits_{j=1}^{v}\Gamma_\x(h_j)\right)\right\rangle\\
		&\subseteq \conv  \left\langle \bigcup\limits_{1\leq i\leq u,1\leq j\leq v}(\Gamma_\x(g_i) +\Gamma_{\x}(h_j)) \right\rangle\\
		&= \conv  \left\langle \bigcup\limits_{1\leq i\leq u,1\leq j\leq v}(\Gamma_\x(g_ih_j)) \right\rangle\\
		&= \Gamma_{\x}(IJ).
	\end{align*}
The assertion is proved.
\end{proof}

As an immediate consequence of \Cref{thm.product}, we obtain the following statement for powers of an ideal.

\begin{corollary} \label{thm.NPPower}
Let $I\subseteq R$ be an ideal. Then, for any $n \in \NN$, we have
$$\Gamma_\x(I^n) = n\Gamma_\x(I).$$
\end{corollary}



We now would like to give criteria that characterize ideals having monomial integral closure via their Newton polyhedra that generalize results in \cite{BAFS02}. Similar to \cite{BAFS02}, if $I$ is an ideal in $R$, then we denote by $I_0$ the ideal generated by those monomials $\x^\a$ with $\a \in \Gamma_\x(I) \cap \ZZ_{\ge 0}^p$. Furthermore, let $K_I$ be the ideal generated by monomials in $R$ that are also in the integral closure $\overline{I}$ of $I$. Note that while it is not hard to show that $I \subseteq I_0$ in a regular local ring (and then follow the argument given in \cite{BAFS02} to show their results in a regular local ring), it is not clear if $I \subseteq I_0$ is true in a regular domain. For our purpose, we only need the following.

\begin{lemma}
	Let $I\subseteq R$ be an ideal. Then, $\Gamma_\x(\overline{I})\subseteq\Gamma_{\x}(\overline{I_0})$. 
\end{lemma}

\begin{proof}
	Let $f$ be an arbitrary element in $\overline{I}$. It is clear that $\supp(f)\in \supp(\overline{I})$. Since $f\in \overline{I}$, then $f$ must satisfy an equation 
	\begin{equation*}
		f^n+a_{1}f^{n-1}+\cdots+a_{n-1}f+a_n=0,
	\end{equation*}
	where $a_i\in I^{i}$ for all $1\leq i\leq n$. For each $i$, there exists an element $h_i\not\in (\x)=(x_1,\ldots,x_p)$ such that $h_ia_i\in (I_0)^i$. Let $h=\prod_{i=1}^{n}h_i$, then we have 
	\begin{equation*}
		(hf)^n+a_{1}h(hf)^{n-1}+\cdots+a_{n-1}h^{n-1}(hf)+a_nh^n=0,
	\end{equation*}
	where $a_ih^i\in (I_0)^i$ for all $i=1,\ldots,n$. It implies that $hf\in \overline{I_0}$. Since $h\not\in(\x)=(x_1,\ldots,x_p)$, then $\supp(f)=\supp(hf)\in \supp(\overline{I_0})$. Therefore, $\supp (\overline{I})\subseteq \supp (\overline{I_0})$, and hence $\Gamma_\x(\overline{I})\subseteq\Gamma_{\x}(\overline{I_0})$. 
\end{proof}

\begin{question}
	\label{ques.IvsI0}
	In a regular domain, is it true that $I \subseteq I_0$?
\end{question}

The following statement exhibits a similar property as that of the Newton polyhedron of monomial ideals in polynomial or formal power series rings.

\begin{lemma}[See {\cite[Lemma 2.1]{BAFS02}}] \label{lem.GammaI0}
	Let $I\subset R$ be an ideal. Then, $\Gamma_\x(I)=\Gamma_\x(\overline{I}) = \Gamma_\x(I_0)$.
\end{lemma}

\begin{proof}
	Observe that, by \Cref{cor.Gamma=NP} and \cite[Theorem 2.3]{HS08}, we have
	$$\Gamma_\x(I_0) = \NP(I_0) = \NP(\overline{I_0}) = \Gamma_\x(\overline{I_0}).$$
	Moreover, it follows from the definition that if $I \subseteq J$ are ideals in $R$, then $\Gamma_\x(I) \subseteq \Gamma_\x(J)$. Therefore,
	$$\Gamma_\x(I) \subseteq \Gamma_\x(\overline{I}) \subseteq \Gamma_\x(\overline{I_0}) = \Gamma_\x({I_0}) = \Gamma_\x(I).$$
	Hence, $\Gamma_\x(I) = \Gamma_\x(\overline{I}) = \Gamma_\x(I_0).$
\end{proof}

Our next theorem generalizes \cite[Theorem 2.3]{BAFS02} to local domains.

\begin{theorem}[Compare to {\cite[Theorem 2.3]{BAFS02}}]
	\label{thm.NNDregular}
	Let $I\subseteq R$ be an ideal. TFAE:
	\begin{enumerate}
		\item $\overline{I}$ is a monomial ideal;
		\item $I_0 \subseteq \overline{I}$;
		\item $\overline{I}=\{f\in R \mid \Gamma_\x(f)\subseteq\Gamma_\x(I)\}$;
		\item $\Gamma_\x(I)=\Gamma_\x(K_I)$.
	\end{enumerate}
\end{theorem}

\begin{proof} Suppose that $\overline{I}$ is a monomial ideal. By \Cref{cor.Gamma=NP} and \Cref{lem.GammaI0}, we have $\NP(\overline{I}) = \Gamma_\x(\overline{I}) = \Gamma_\x(I_0) = \NP(I_0)$. It follows from \cite[Theorem 2.3]{HS08} that $I_0 \subseteq \overline{I_0} = \overline{I}$. Thus, $I_0 \subseteq \overline{I}$. Therefore, $(1) \Rightarrow (2)$.

To see $(2) \Rightarrow (3)$, let $f$ be any element in $\overline{I}$. Then we have $\Gamma_\x(f)\subseteq \Gamma_\x(\overline{I})=\Gamma_\x(I)$, hence $\overline{I}\subseteq\{f\in R \mid \Gamma_\x(f)\subseteq\Gamma_\x(I)\}$. To obtain the reverse inclusion, consider an arbitrary element $f\in R$ with $\Gamma_\x(f)\subseteq \Gamma_\x(I)$. Then, all monomials in a minimal monomial representation of $f$ are in $I_0$ by definition. Therefore, $f\in I_0 \subseteq \overline{I}$. That is,
$$\{f\in R \mid \Gamma_\x(f)\subseteq\Gamma_\x(I)\} \subseteq\overline{I}.$$

We proceed to prove $(3)\Rightarrow (1)$. Consider any element $g \in \overline{I}$. Clearly, $\a \in \supp(g)$ implies that $\Gamma_\x(\x^\a) \subseteq \Gamma_\x(g) \subseteq \Gamma_\x(\overline{I}) = \Gamma_\x(I)$. Thus, for any $\a \in \supp(g)$,
$$\a \in \{f \in R \mid \Gamma_\x(f) \subseteq \Gamma_x(I)\}.$$
It follows that $\a \in \overline{I}$, by the hypothesis. Particularly, by running $g$ through a system of generators of $\overline{I}$, we get a monomial generating set for $\overline{I}$; that is, $\overline{I}$ is a monomial ideal.

Finally, we show that $(1) \Leftrightarrow (4)$. The implication $(1) \Rightarrow (4)$ is obvious as $\Gamma_\x(I)=\Gamma_\x(\overline{I})$, by \Cref{lem.GammaI0}, and $(1)$ implies that $K_I = \overline{I}$. To exhibit $(4) \Rightarrow (1)$, consider $g \in \overline{I}$ and any $\a \in \supp(g)$. By \Cref{lem.GammaI0} and the hypothesis, we have $\a \in \Gamma_\x(\overline{I}) = \Gamma_\x(I) = \Gamma_\x(K_I)$. It then follows from \cite[Theorem 2.3]{HS08} that $\x^\a \in \overline{K_I}\subseteq \overline{\overline{I}} = \overline{I}$. Again, run $g$ through a collection of generators of $\overline{I}$, we obtain a monomial generating set of $\overline{I}$, showing that $\overline{I}$ is a monomial ideal. The proof is complete.
\end{proof}


\begin{corollary}\label{cor.IntegralClosureofNNDidealsPoly}
	Let $I,J$ be ideals of $R$ such that their integral closures are monomial ideals. TFAE:
	\begin{enumerate}
		\item $\overline{I}=\overline{J}$;
		\item $\Gamma_\x(I)=\Gamma_\x(J)$.
	\end{enumerate}
\end{corollary}

\begin{proof}
	(1) $\Rightarrow$ (2) is obvious. We will show that (2) $\Rightarrow$ (1). Let $f$ be an element in $\overline{I}$. Then by (3) of Theorem \ref{thm.NNDregular}, $\Gamma_\x(f)\subseteq \Gamma_\x(I)$. Since $\Gamma_\x(I)=\Gamma_\x(J)$, then $\Gamma_\x(f)\subseteq \Gamma_\x(J)$. Using part (3) of Theorem \ref{thm.NNDregular} again, we have $f\in \overline{J}$. Thus $\overline{I}\subseteq\overline{J}$. Similarly, $\overline{J}\subseteq \overline{I}$. Therefore $\overline{J}= \overline{I}$.
\end{proof}


\section{Newton non-degenerate ideals in regular local rings} \label{sec.NND}

In this section, we introduce the concept of \emph{Newton non-degenerate} (NND) ideals in regular local rings and demonstrate that this notion can be characterized by properties of the integral closures of the ideals. Throughout this section, $(R,\mm)$ denotes a regular local ring of dimension $d$ which contains a field. We also fix a \emph{regular system of parameters} (r.s.o.p) $\x = x_1, \dots, x_d$ of $R$. As remarked in \cite{HS08}, $\x$ is also a generalized regular system of parameters in $R$, and so the theory of monomials, monomial ideals and Newton polyhedra with respect to $\x$ is the same as that in \Cref{sec.Monomials}.

\medskip

\noindent \textbf{Notation.} Let $\Delta$ be a close convex set in $\RR_{\ge 0}^d$. Let $C(\Delta)$ be the infinite cone over $\Delta$ with the origin as its vertex; i.e., $C(\Delta)$ consists of half-rays emanating from the origin and going through points in $\Delta$. Set
$$R_\Delta = \{g \in R ~\big|~ \supp(g) \subseteq C(\Delta) \cap \ZZ_{\ge 0}^d\}.$$
Observe that, since $\Delta$ is convex, $C(\Delta) \cap \ZZ_{\ge 0}^d$ is a subsemigroup of $\ZZ_{\ge 0}^d$. Thus, $R_\Delta$ is a subring of $R$. Moreover, $R_\Delta$ is a local ring with a unique maximal ideal
$$\mm_\Delta = \{g \in R_\Delta ~\big|~ \0 = (0, \dots, 0) \not\in \supp(g)\}.$$
For an element $g \in R$ with its unique monomial representation with respect to $\x$ being
$g = \sum\limits_{\a \in \supp(g)} c_\a\x^\a,$
set $$g_\Delta = \sum_{\a \in \supp(g) \cap C(\Delta)} c_\a \x^\a.$$

\begin{definition} \label{def.NND}
An ideal $I \subseteq R$ is called \emph{Newton non-degenerate} (w.r.t. $\x$) if there exists a system of generators $g_1, \dots, g_s$ of $I$ such that, for each compact face $\Delta \subseteq \Gamma_\x(I)$, the ideal $I_\Delta$ generated by $g_{1\Delta}, \dots, g_{s\Delta}$ is an $\mm_\Delta$-primary ideal in $R_\Delta$.
\end{definition}

\begin{remark} \label{rmk.K}
In the case when $R = \CC[[x_1, \dots, x_d]]$ is the ring of formal power series over the complex number, \cite[Theorem 6.2]{Kou76} shows that the ideal $I_\Delta$ is $\mm_\Delta$-primary if and only if for each compact face $\Delta'$ of $\Delta$, the system of equations $g_{1\Delta'} = \dots = g_{s\Delta'} = 0$ has no common solution in $(\CC \setminus 0)^d$. This justifies the terminology ``Newton non-degenerate''.
\end{remark}

The main result of this section is stated as follows.

\begin{theorem} \label{thm.NNDequiv}
Let $(R,\mm)$ be a regular local ring of dimension $d$ with a regular system of parameters $\x = x_1, \dots, x_d$. Let $I \subseteq R$ be an ideal. Then, $I$ is NND if and only if its integral closure $\bar{I}$ is a monomial ideal in $\x$.
\end{theorem}

To prove \Cref{thm.NNDequiv}, we will need the following lemma.

\begin{lemma} \label{lem.complete}
Let $\widehat{R}$ be the completion of $R$ and let $\widehat{\mm} = \mm \widehat{R}$ be its maximal ideal. Let $I \subseteq R$ be an ideal. Then, $I\widehat{R}$ is $\widehat{\mm}$-primary if and only if $I$ is $\mm$-primary.
\end{lemma}

\begin{proof}
	Suppose that $I\widehat{R}$ is $\widehat{\mm}$-primary. We just need to show that $\mm\subseteq \sqrt{I}$. Let $x\in \mm$ be an arbitrary element. Then, $x\in \widehat{\mm}=\sqrt{I\widehat{R}}$, i.e., $x^k\in I\widehat{R}$ for some $k\in \NN$. Thus, $x^k\in  I\widehat{R}\cap R=I$, since the map $R \rightarrow \widehat{R}$ is faithfully flat. That is, $x\in \sqrt{I}$. Therefore, $\sqrt{I}=\mm$.
	
	Conversely, suppose that $\sqrt{I}=\mm$. Then, we have $\widehat{\mm}=\mm\widehat{R}=\sqrt{I}\widehat{R}\subseteq \sqrt{I\widehat{R}}$. It follows that $I\widehat{R}$ is $\widehat{\mm}$-primary.
\end{proof}


\begin{corollary}\label{cor_m_primary_between_R_and_hatR}
	Let $I$ be an ideal of $R$. Then, for each compact face $\Delta\subseteq\Gamma_\x(I)$, $I(\widehat{R})_{\Delta}$ is $\widehat{\mm}_{\Delta}$-primary in $(\widehat{R})_{\Delta}$ if and only if $I_{\Delta}$ is $\mm_{\Delta}$-primary in $R_{\Delta}$.
\end{corollary}

\begin{proof}
	It suffices to show that for any ideal $I$ of $R_{\Delta}$, $I(\widehat{R})_{\Delta}\cap R_{\Delta}=I$. Indeed, we have
	\begin{equation*}
	I \subseteq	I(\widehat{R})_{\Delta}\cap R_{\Delta}\subseteq I\widehat{R}\cap R_{\Delta}=I\widehat{R}\cap (R\cap R_{\Delta})=(I\widehat{R}\cap R)\cap R_{\Delta}= (IR)\cap R_{\Delta} =I. \qedhere
	\end{equation*}
\end{proof}

We are now ready to prove \Cref{thm.NNDequiv}.

\begin{proof}[{\bf Proof of \Cref{thm.NNDequiv}.}] Note that $R$ is a regular local ring, so $\widehat{R} \simeq k[[X_1, \dots, X_d]]$, a formal power series ring, in which $k \simeq R/\mm$. Clearly, $\mathbf{X} = X_1,\dots, X_d$ form a regular system of parameters in $\widehat{R}$.

Suppose first that $I$ is an NND ideal in $R$. By Corollary \ref{cor_m_primary_between_R_and_hatR}, $I\widehat{R}$ is an NND ideal in $\widehat{R}$. Thus, it follows from \cite[Corollary 2.1]{BAFS02} that the integral closure $\overline{I\widehat{R}}$ is generated by monomials; note that \cite[Corollary 2.1]{BAFS02} was stated for ideals of finite codimension, but the the same proof goes through for any ideal. Since $R \rightarrow \widehat{R}$ is faithfully flat, we have $\overline{I}=\overline{I\widehat{R}}\cap R$. Hence, $\overline{I}$ is also a monomial ideal in $x_1,\dots,x_d$.

Conversely, suppose that $\overline{I}$ is a monomial ideal. Again, since $R \rightarrow \widehat{R}$ is faithfully flat, $\overline{I\widehat{R}} = \overline{I} \widehat{R}$ is a monomial ideal. By the same proof of \cite[Corollary 2.1]{BAFS02}, we have that $I\widehat{R}$ is an NND ideal in $\widehat{R}$. That is, for each compact face $\Delta\subseteq\Gamma_{\mathbf{X}}(I\widehat{R})$, the ideal $(I\widehat{R})_{\Delta}$ is $\widehat{\mm}_{\Delta}$-primary in $\widehat{R}_{\Delta}$. Moreover, it follows from \Cref{thm.NPgenerators} that
$$\Gamma_{\mathbf{X}}(I\widehat{R}) = \conv\left\langle \Gamma_{\mathbf{X}}(g_1) \cup \dots \Gamma_{\mathbf{X}}(g_s) \right\rangle = \conv \left\langle \Gamma_\x(g_1) \cup \dots \cup \Gamma_\x(g_s)\right\rangle = \Gamma_\x(I),$$
where $g_1, \dots, g_s$ is a generating set of $I$.
Therefore, Corollary \ref{cor_m_primary_between_R_and_hatR} now implies that $I_{\Delta}$ is $\mm_{\Delta}$-primary in $R_{\Delta}$ for any compact face $\Delta \subseteq \Gamma_\x(I)$. Hence, $I$ is an NND ideal in $R$.
\end{proof}

The following remark shows that if one attempts to define a notion of NND ideals in the nonlocal ring, for instance, using the condition on non-common solutions in $(\CC \setminus 0)^d$ as in Remark \ref{rmk.K}, they might lose the property of having monomial integral closure. Furthermore, this remark also shows that having a monomial integral closure in the completion ring does not imply having a monomial integral closure in the original ring.

\begin{remark}\label{rmk.NonLocal}
Consider the ideal
$$I = (x^5+xy^3, y^5+x^3y) \subseteq R = \CC[x,y].$$
It can be seen that $I$ satisfies the condition that the system of equations $g_{1\Delta'} = g_{2\Delta'} = 0$ has no common solution in $(\CC \setminus 0)^2$ with $g_1=x^5+xy^3$, $g_2=y^5+x^3y$ as in Remark \ref{rmk.K}. One observes that the integral closure of $IS$, where $S = \CC[[x,y]]$, is a monomial ideal in $S$, so $IS$ is NND in $S$ and satisfies the condition on having no common solutions above by Theorem \ref{thm.NNDequiv}. However, the integral closure of $I$ in $R$ is
$$\overline{I} = (y^5+x^3y, x^5+xy^3, x^2y^4-xy^3, x^3y^3-x^2y^2, x^4y^2-x^3y)$$
which is not a monomial ideal.
\end{remark}

\Cref{thm.NNDequiv} provides a simple method for checking whether an ideal is NND. In particular, it gives the following properties of NND ideals.

\begin{corollary} \label{NNDSumInt}
Let $I$ and $J$ be NND ideals in $R$.
\begin{enumerate}
    \item The ideal $IJ$ is NND. 
    \item If, in addition, $\text{ht}(I+J) = \ell(I) + \ell(J)$, where $\ell(I)$ and $\ell(J)$ are the analytic spreads of $I$ and $J$, respectively, then $I \cap J$ is an NND ideal.
\end{enumerate}
\end{corollary}

\begin{proof} (1) The argument is similar to \cite[Corollary 2.4]{BAFS02}. To see (2), observe that $R$ is a Cohen Macaulay ring. Thus, $R_{\pp}$ is also Cohen-Macaulay for any prime ideal $\pp \subseteq R$. It implies that the completion $\widehat{R_{\pp}}$ is Cohen-Macaulay, namely  $\widehat{R_{\pp}}$ is equidimensional. It follows that $R$ is formally equidimensional. Therefore, by \cite[Exercise 10.24]{SH06}, we have $\overline{IJ}=\overline{I}\cap \overline{J}$. Since
	\begin{equation*}
		 \overline{IJ}\subseteq \overline{I\cap J} \subseteq\overline{I}\cap \overline{J},
	\end{equation*}
	then we must have $\overline{I\cap J}=\overline{IJ}=\overline{I}\cap \overline{J}$, namely it is a monomial ideal. Thus $I\cap J$ is NND.
\end{proof}

\begin{remark}
    \label{rmk.NNDSum}
        Although verifying that the sum of two NND ideals $I$ and $J$ in a same ring $R$ is NND is not an easy task, it can still be proven in some special cases. Consider the following case: Let $(S,\mathfrak{m}),(T,\mathfrak{n})$ be regular local domains with the same residue field $k$. Let $I \subseteq R$ and $J \subseteq S$ be NND ideals, we can show that $I+J$ is an NND ideal in $R\otimes_kS$ by using the Rees packages of $I$ and $J$ as in \cite[Theorem 4.8]{SSTJ24}.
\end{remark}

We conclude this section with the following result demonstrating that if a power $I^n$, for some $n \in \NN$, of an ideal $I$ is NND then so is $I$.

\begin{theorem}\label{thm.power.is.NND}
	Let $I \subseteq R$ be an ideal. If $I^n$ is NND, for some $n>0$, then $I$ is also NND.
\end{theorem}

\begin{proof}
	Since $I^n$ is NND, we have
$$e(I^n)=e((I^n)_0).$$
Moreover, since $(I_0)^n\subseteq (I^n)_0$, we also have
$$e((I^n)_0) \leq e(I_0)^n=n^d e(I_0).$$
It follows that
	\begin{equation*}
		n^de(I)=e((I^n)_0)\leq e(I_0)^n=n^d e(I_0).
	\end{equation*}
That is, $e(I)\leq e(I_0)$.

On the other hand, since $I\subseteq I_0$, we have $e(I)\geq e(I_0)$. Thus, $e(I)=e(I_0)$. This implies that $I$ and $I_0$ have the same integral closure. In particular, $I$ must be NND is asserted.
\end{proof}


\section{Multiplicity and graded families of Newton non-degenerate ideals} \label{sec.Multiplicity}

In this section, we examine graded families of ideals in regular local rings and the ``Multiplicity $=$ Volume'' formula. We will show that another type of convex bodies can be used in place of the Newton-Okounkov bodies in this formula if the family contains a subfamily of Newton non-degenerate ideals. As in \Cref{sec.NND}, throughout this section, $(R,\mm)$ is a regular local ring of dimension $d$, which contains a field, and $\x = x_1, \dots, x_d$ denotes a regular system of parameters in $R$.

We shall briefly recall the notion of graded families of ideals and their associated Rees algebras.

\begin{definition} \label{def.gradedFamily} \quad
\begin{enumerate}
     \item A collection $\I = \{I_n\}_{n \in \NN}$ of ideals in $R$ is called a \emph{graded} family if $I_p \cdot I_q \subseteq I_{p+q}$ for any $p, q \in \NN$. A \emph{filtration} is a graded family in which $I_p \supseteq I_{p+1}$ for any $p \in \NN$.
     \item If $\I = \{I_n\}_{n \in \NN}$ is a graded family of ideals in $R$, then the \emph{Rees algebra} of $\I$ is defined to be
         $$\R(\I) = \bigoplus_{n \ge 0} I_n t^n \subseteq R[t],$$
         where, by convention, $I_0 = R$.
         The family $\I$ is said to be \emph{Noetherian} if its Rees algebra is a Noetherian ring.
\end{enumerate}
\end{definition}

A graded family $\I$ is called \emph{Noetherian} if its Rees algebra $\R(\I)$ is a Noetherian ring. The property of having Noetherian Rees algebras is fundamental, as it enables many powerful results to hold and is closely tied to Nagata’s work \cite{Nag59} on Hilbert’s 14th problem. Numerous examples of non-Noetherian \emph{symbolic} Rees algebras have been studied in the literature (cf. \cite{Cut91, Hun82, Rob85}). For graded families of monomial ideals, the Noetherian property of their Rees algebras was characterized in \cite{HN24}. This characterization was in terms of the limiting bodies associated to these families of ideals.

We shall extend the construction of limiting bodies and results of \cite{HN24} to provide a classification of graded families of Newton non-degenerate ideals in a regular local ring whose Rees algebras are Noetherian.

\begin{definition} \label{def.limiting}
Let $\I = \{I_n\}_{n \in \NN}$ be a graded family of ideals in $R$. The \emph{limiting body} of $\I$ is defined as
$$\C(\I) = \bigcup_{n \in \NN} \frac{1}{n} \Gamma_\x(I_n) \subseteq \RR_{\ge 0}^d.$$
\end{definition}

The following lemma shows that $\C(\I)$ is a convex set. However, it is not necessarily closed.

\begin{lemma} \label{lem.CIconvex}
	The limiting body $\C(\I)$ of a graded family $\I = \{I_n\}_{n \in \NN}$ of ideals in $R$ is a convex set.
\end{lemma}

\begin{proof}
	We have
	$(I_n)^k \subseteq I_{nk} $ for all $k \in \mathbb{N}$. This implies that, for all $k \in \NN$,
	\begin{equation*}
		\frac{1}{n}\Gamma_\x(I_n)\subseteq\frac{1}{nk}\Gamma_\x(I_{nk}).
	\end{equation*}
Consider arbitrary $u, v \in \C(\I)$. Suppose that $u \in \frac{1}{a}\Gamma_\x(I_a)$ and $v \in
	\frac{1}{b}\Gamma_\x(I_b)$ for some $a,b \in \NN$. Then, $u, v \in \frac{1}{ab}\Gamma_\x(I_{ab})$. It follows, since $\Gamma_\x(I_{ab})$ is a convex set, that any convex combination of $u,v$ is also in $\frac{1}{ab}\Gamma_\x(I_{ab})$ and, hence, in $\C(\I)$. Therefore, $\C(\I)$ is a convex set.
\end{proof}

The following results were proved in \cite{HN24} for graded families of monomial ideals, but the proofs carry verbatim for families of arbitrary ideals.

\begin{lemma}[{\cite[Theorem 3.1]{HN24}}]
	\label{thm.polyhedronPoly}
	Let $\I=\{I_n\}_{n\in \NN}$ be a graded family of ideals in $R$. TFAE:\begin{enumerate}
		\item The limiting body $\C(\I)$ is a polyhedron;
		\item There exists an integer $c$ such that $\overline{\C(\I)}=\frac{1}{c}\Gamma_\x(I_c)$;
		\item There exists an integer $c$ such that
		$\frac{1}{c}\Gamma_\x(I_c)=\frac{1}{kc}\Gamma_\x(I_{kc})$ for all $k\in \mathbb{N}$.
	\end{enumerate}
\end{lemma}

\begin{proof}
	The proof goes in exactly the same line of arguments as that of \cite[Theorem 3.1]{HN24}.
\end{proof}

\begin{lemma}[{\cite[Theorem 3.4]{HN24}}]
	\label{thm.Noetherian}
	Let $\I=\{I_n\}_{n\in \NN}$ be a graded family of ideals in $R$ and let $\overline{\I}=\{\overline{I_n}\}_{n\in \NN}$. TFAE:
	\begin{enumerate}
		\item There exists an integer $c$ such that $\overline{I_c^k}=\overline{I_{kc}}$ for all $k\in \mathbb{N}$;
		\item $\mathcal{R}(\overline{\I})$ is Noetherian;
		\item $\mathcal{R}(\I)$ is Noetherian.
	\end{enumerate}
\end{lemma}

\begin{proof}
	The proof goes in exactly the same line of arguments as that of \cite[Theorem 3.4]{HN24}.
\end{proof}

The classification for Noetherian property of the Rees algebra of a graded family of NND ideals in a regular local ring is given as follows.

\begin{theorem}
	\label{thm.Noe.property.NND.ideal}
	Let $\I=\{I_n\}_{n\in \NN}$ be a graded family of NND ideals in $R$ and let $\overline{\I}=\{\overline{I_n}\}_{n\in \NN}$. TFAE:
	\begin{enumerate}
		\item $\I$ satisfies any of the equivalent conditions in \Cref{thm.polyhedronPoly};
		\item There exists an integer $c$ such that $\overline{I_c^k}=\overline{I_{kc}}$ for all $k\in \mathbb{N}$;
		\item $\mathcal{R}(\overline{\I})$ is Noetherian;
		\item $\mathcal{R}(\I)$ is Noetherian.
	\end{enumerate}
\end{theorem}

\begin{proof}
The equivalence between (2),(3), and (4) was the content of \Cref{thm.Noetherian}. It remains to show that (1) and (2) are equivalent.

Indeed, observe that for any $k,c \in \NN$, $I_c^k$ is an NND ideal by \Cref{NNDSumInt}. Thus, by \Cref{thm.NNDequiv}, both $\overline{I_c^k}$ and $\overline{I_{kc}}$ are monomial ideals. It now follows from \Cref{cor.IntegralClosureofNNDidealsPoly} that $\overline{I_c^k} = \overline{I_{kc}}$ if and only if $\Gamma_\x(I_c^k) = \Gamma_\x(I_{kc})$. Hence, (2) occurs if and only if there exists $c \in \NN$ such that $\Gamma_\x(I_c^k) = \Gamma_\x(I_{kc})$ for all $k \in \NN$, which is one of the equivalent conditions in \Cref{thm.polyhedronPoly}.
\end{proof}

\begin{example}
    Consider the family $\mathcal{I}=\{I_k\}_{k\in \mathrm{N}}$ in $k[[x,y]]$ with 
	\begin{equation*}
		I_k=(x^{\lceil k/2\rceil+1}+y^{\lceil k/2\rceil+1}, x^iy^j), 
	\end{equation*}
	where $i,j\in \mathbb{N}\setminus\{0\}$ and $i+j=\lceil k/2\rceil+1$. We claim that
    \begin{enumerate}
        \item $\mathcal{I}=\{I_k\}_{k\in \mathrm{N}}$ is a graded family of NND ideals.
        \item $\mathcal{R}(\I)$ is not Noetherian.
    \end{enumerate}
    \begin{proof}
        For the first statement, we show the following stronger statement that for any fixed integer $m$, the ideal $J_m=(x^m+y^m, x^iy^j) \subset k[[x,y]]$ such that $i,j\in \mathbb{N}\setminus\{0\}, i+j=m$ is an NND ideal. Indeed, since 
	\begin{equation*}
		(x^m)^m-\left[x^m+y^m+\sum_{j=1}^{m-2}x^jy^{m-j}\right](x^m)^{m-1}+\sum_{j=1}^{m-2}(x^{m-1}y)^{m-j}(x^m)^j+(x^{m-1}y)^m=0,
	\end{equation*}
	then $x^m\in \overline{J_m}$. It also implies that $y^m\in \overline{J_m}$, and hence $\overline{J_m}=(x,y)^m$. 
    
    \vspace{0.5em}
    For the second statement, note that 
	\begin{equation*}
		\dfrac{1}{k}\Gamma_+(I_k)=\mathrm{conv}\left\{\left(\dfrac{\lceil k/2\rceil+1}{k},0\right),\left(0,\dfrac{\lceil k/2\rceil+1}{k}\right)\right\}+\mathbb{R}^2_{\geq 0}.
	\end{equation*}
	Since $\lim\limits_{n\to \infty}\dfrac{\lceil k/2\rceil+1}{k}=\dfrac{1}{2}$, $\mathcal{C}(\mathcal{I})$ is the interior of $\conv\left\{ (\frac{1}{2},0),(0,\frac{1}{2})\right\}+\mathbb{R}^2_{\geq 0}$ together with two open rays $(\frac{1}{2},\infty)\times 0$ and $0\times (\frac{1}{2},\infty)$. Because $\mathcal{C}(\mathcal{I})$ is not a polyhedron, by Theorem \ref{thm.Noe.property.NND.ideal}, $\mathcal{R}(\I)$ is not Noetherian. 
    \end{proof} 
\end{example}

\begin{theorem}\label{thm.Mult.mono}
Let $(R,\mm)$ be a regular local ring of dimension $d$ and let $\x = x_1, \dots, x_d$ be a regular system of parameters in $R$. If $I$ is an $\mm$-primary monomial ideal in $\x$, then $e(I)=d!\text{co-vol}_d (\Gamma_{\x}(I))$.    
\end{theorem}

\begin{proof}
    Note that $\widehat{R}\cong k[[X_1,\ldots,X_d]]$ where $k\cong R/\mm$. Since $I$ is an $\mm$-primary monomial ideal in $R$, then $\widehat{I}=I\widehat{R}$ is an $\widehat{\mm}$-primary ideal in $\widehat{R}$. It follows that 
    \begin{equation*}
        e(I)=e(\widehat{I})=d!\text{co-vol}_d (\text{NP}(\widehat{I})),
    \end{equation*}
    where the first equality is because the length $\ell_R(R/J)$ is unaffected by completion for any $\mm$-primary ideal $J$, and the second equality follows from \cite[p. 131]{Tei88}. Moreover, by construction, $\Gamma_{\x}(I)$ and $\text{NP}(\widehat{I})$ coincide. Therefore, $e(I)=d!\text{co-vol}_d (\Gamma_{\x}(I))$ as desired.  
\end{proof}

We are now ready to state and prove the last main result of the paper.

\begin{theorem} \label{thm.Mult}
Let $(R,\mm)$ be a regular local ring of dimension $d$ and let $\x = x_1, \dots, x_d$ be a regular system of parameters in $R$. Let $\I = \{I_n\}_{n \in \NN}$ be a Noetherian graded family of $\mm$-primary ideals in $R$. Let $c$ be an integer such that $\overline{I_c^k} = \overline{I_{kc}}$, for all $k \in \NN$, as in \Cref{thm.Noetherian}. Then, TFAE:
\begin{enumerate}
\item $e(\I) = d! \text{co-vol}_d (\C(\I))$, and
\item $I_{kc}$ is an NND ideal for every $k \in \NN$.
\end{enumerate}
\end{theorem}

\begin{proof}
	We first prove (1) $\Rightarrow$ (2). Suppose that (1) holds. Recall that, for an ideal $J$, the ideal generated by monomials $\{\x^\a ~\big|~ \a \in \Gamma_\x(J) \cap \ZZ_{\ge 0}^d\}$ is denoted by $J_0$. By \Cref{cor.Gamma=NP} and \Cref{thm.Mult.mono}, we have  
	\begin{align*}
		e(\I)&= d!\text{co-vol}_d\left(\C(\I)\right)
		= d!\text{co-vol}_d\left(\frac{1}{c}\Gamma_\x(I_c)\right)\\
		&=d!\frac{1}{c^d}\text{co-vol}_d\left(\Gamma_\x(I_c)\right)
		=d!\frac{1}{c^d}\text{co-vol}_d\left(\Gamma_\x((I_c)_0)\right)\\
		&=\frac{1}{c^d} d! \text{co-vol}_d(\NP((I_c)_0)) = \frac{e((I_c)_0)}{c^d}.
	\end{align*}
	
On the other hand, \cite[Theorem 6.5]{Cut13} gives
$$e(\I)=\lim\limits_{n\to \infty}\frac{e(I_n)}{n^d}.$$
It follows that
	\begin{equation} \label{eq.eI}
		e(\I)=\lim\limits_{k\to \infty}\frac{e(I_{kc})}{(kc)^d}=\lim\limits_{k\to \infty}\frac{e(\overline{I_{kc}})}{(kc)^d}=\lim\limits_{k\to \infty}\frac{e(\overline{I_c^k})}{(kc)^d}=\lim\limits_{k\to \infty}\frac{e(I_c^k)}{(kc)^d}=\lim\limits_{k\to \infty}\frac{k^de(I_c)}{k^dc^d}=\frac{e(I_c)}{c^d}.
	\end{equation}
Thus,
$$e((I_c)_0)=e(I_c).$$
Since $I_c\subseteq (I_c)_0$, this implies that $I_c$ is a reduction of $(I_c)_0$; that is, $\overline{I_c} = \overline{(I_c)_0}$. It then follows from \cite[Theorem 2.3]{HS08} that $\overline{I_c}$ is generated by monomials.

Observe further that, by \Cref{thm.NPPower}, for any $k \in \NN$, we have
$$\Gamma_\x((I_c)_0^k) = k\Gamma_\x((I_c)_0) = k \Gamma_\x(I_c) = \Gamma_\x(I_c^k) = \Gamma_\x(\overline{I_c^k}) = \Gamma_\x(\overline{I_{kc}}) = \Gamma_\x(I_{kc}) = \Gamma_\x((I_{kc})_0).$$
Together with \Cref{cor.Gamma=NP}, since $(I_c)_0^k$ and $(I_{kc})_0$ are both monomial ideals in $R$, this implies that $\NP((I_c)_0^k) = \NP((I_{kc})_0)$; that is, $\overline{(I_c)_0^k} = \overline{(I_{kc})_0}$. Particularly, for any $k \in \NN$, we get
$$e(I_{kc}) = e(I_c^k) = k^de(I_c) = k^de((I_c)_0) = e((I_c)_0^k) = e((I_{kc})_0).$$
Since $I_{kc} \subseteq (I_{kc})_0$, by a similar reasoning as that for $I_c$, we conclude that $\overline{I_{kc}} = \overline{(I_{kc})_0}$ and $\overline{I_{kc}}$ is a monomial ideal for any $k \in \NN$. Hence, $I_{kc}$ is an NND ideal for any $k \in \NN$, by \Cref{thm.NNDequiv}.

We finally establish (2) $\Rightarrow$ (1). Suppose $I_{kc}$ is NND for all $k \in \NN$. By \Cref{thm.NNDequiv}, we have that $\overline{I_c}$ is a monomial ideal. This, together with the proof of \Cref{thm.NNDregular}, implies that $\overline{(I_c)_0} = \overline{I_c}$. Thus,
$$e(I_c)=e((I_c)_0).$$
Therefore, as in (\ref{eq.eI}) and \Cref{thm.Mult.mono}, we get
	\begin{equation*}
		e(\I)=\frac{e(I_c)}{c^d}=\frac{e((I_c)_0)}{c^d}= d!\text{co-vol}_d \left(\frac{1}{c} \NP(I_c)\right) = d!\text{co-vol}_d\left(\frac{1}{c}\Gamma_\x(I_c)\right)=d!\text{co-vol}_c\left(C(\I)\right).
	\end{equation*}
	The proof is complete.
\end{proof}

As a direct consequence of Theorem \ref{thm.Mult} for $\mathcal{I}=\{I^n\}_{n\in \NN}$ we have the following.
\begin{corollary}
   Let $(R,\mm)$ be a regular local ring of dimension $d$ and let $\x = x_1, \dots, x_d$ be a regular system of parameters in $R$. Let $I$ be an $\mm$-primary ideal in $R$. Then, TFAE:
\begin{enumerate}
    \item $e(I)=d!\text{co-vol}_d (\Gamma_{\x}(I))$, and 
    \item $I$ is an NND ideal.
\end{enumerate} 
\end{corollary}

\begin{remark}
    In general, we always have the following inequality 
    \begin{equation*}
        e(\I)\geq e(\I_0)=n!\text{co-vol}_n (\C(\I_0))=n!\text{co-vol}_n (\C(\I)),
    \end{equation*}
    where $\I_0=\{(I_n)_0\}_{n\in \NN}$ and recall that $I_0$ the ideal generated by those monomials $\x^\a$ with $\a \in \Gamma_\x(I) \cap \ZZ_{\ge 0}^p$. It follows from Theorem \ref{thm.Mult} that for a Noetherian graded family $\mathcal{I}=\{I_n\}_{n\in \NN}$, if $I_{kc}$ is not an NND ideal for \textit{some} $k\in \NN$ where $c$ is a constant given in Theorem \ref{thm.Mult}, then so is $I_c$ by Theorem \ref{thm.power.is.NND}, and we have a strict inequality $e(\I)>n!\text{co-vol}_n (\C(\I))$.
\end{remark}

The following example illustrates that when $\I$ does not contain a subfamily of NND ideals, the equality $e(\I) = d! \text{co-vol}_d (\C(\I))$ may not necessarily hold.

\begin{figure}[h]
	\caption{$\Gamma_\x(I)$ for $I = (x^4+y^4, xy^2+x^2y) \subseteq \CC[[x,y]]$.}
	\centering
	\includegraphics[width=0.4\textwidth]{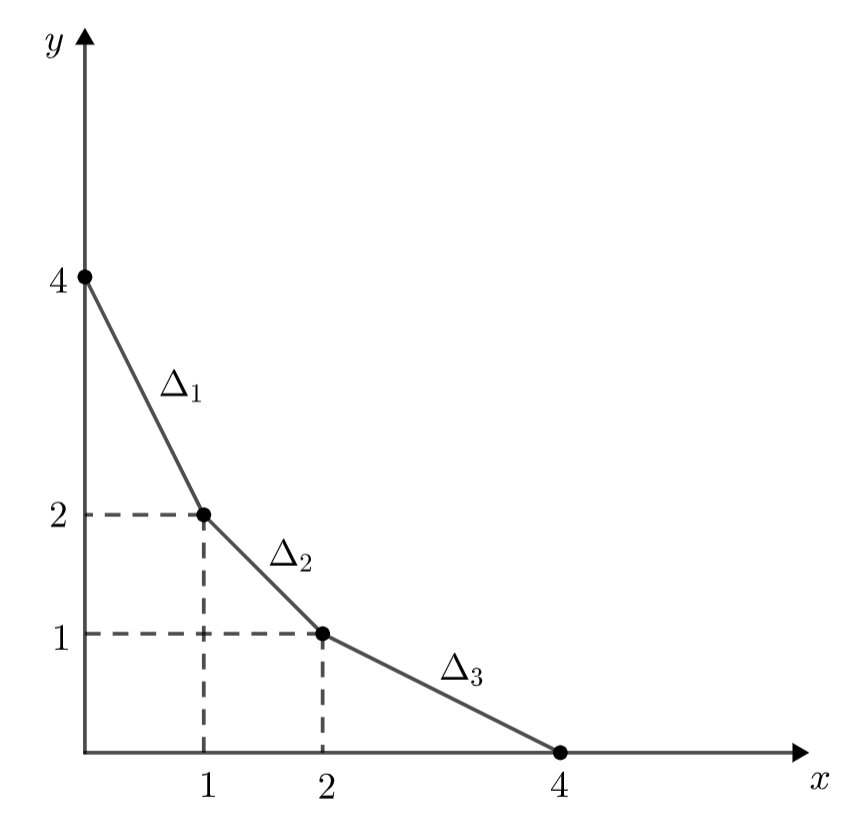}
	\label{non_nnd_ideal}
\end{figure}

\begin{example} \label{ex.NotNNDvolume}
Let $I=(x^4+y^4,xy^2+x^2y)\in R = \CC[[x,y]]$. Consider the $\mm$-primary graded family $\mathcal{I}=\{I^n\}_{n \in \NN}$. It can be seen that
\begin{itemize}
\item $\Gamma_\x(I)$ is given as depicted in Figure \ref{non_nnd_ideal} with three compact faces $\Delta_1, \Delta_2$ and $\Delta_3$;
\item it follows from \Cref{thm.NPPower} that, for any $n \in \NN$, $\Gamma_\x(I^n) = n \Gamma_\x(I)$ also has three compact faces, among which $\Delta_{2n} = \conv\langle \{(n,2n), (2n,n)\}\rangle$.
\end{itemize}
We observe that the system of equations $(x^4+y^4)_{\Delta_2} = (xy^2+x^2y)_{\Delta_{2}} = 0$ has $(t,-t)$, for $t \not= 0$, as common solutions in $(\CC\setminus \{(0,0)\})^2$. This, by \cite[Theorem 6.2]{Kou76} (see also \Cref{rmk.K}), implies that $I$ is not NND. Therefore, $I^n$ is not an NND ideal for any $n \in \NN$, by Theorem \ref{thm.power.is.NND}.

On the other hand, it can be seen that
$$\C(\I) = \conv\left\langle \{(4,0), (2,1), (1,2), (0,4)\}\right\rangle + \RR_{\ge 0}^2.$$
Thus,
$$e(\I) = 12 \not= 11 = 2!\text{co-vol}_2 (\C(\I)).$$
\end{example}

The following example illustrates the use of \Cref{thm.Mult}.

\begin{example} \label{ex.NNDVolume}
Let $R = \CC[[x,y]]$ and consider the following graded family $\I = \{I_n\}_{n \in \NN}$ of $\mm$-primary ideals in $R$:
	\begin{equation*}
			I_1=(x+y),
			I_{2n}=(x^2+y^2,xy)^{n} \text{\space for all \space} n\in \NN,
			I_{2n+1}=I_1I_{2n} \text{\space for all \space} n\in \NN.
	\end{equation*}

Observe that
\begin{itemize}
\item Since $I_2$ is NND, it follows from \Cref{NNDSumInt} that $I_{2n} = I_2^n$ is NND for all $n\geq 1$. On the other hand, for any $n\geq 1$, $I_{2n+1}$, is not NND. This is because over the (only) compact face $\mathrm{conv}\{(2n+1,0),(0,2n+1)\}$, the restricted system of equations will be reduced to the systems with equations in the form $(x+y)(x^2+y^2)^k(xy)^{n-k}=0$ that have common solutions $(t,-t)$ for $t\neq 0$, and so \cite[Theorem 6.2]{Kou76} (see also \Cref{rmk.K}) applies.
\item Since $I_{2n} = I_2^n$ for all $n \in \NN$, the second Veronese subalgebra $\R^{[2]}(\I)$ of $\R(\I)$ is Noetherian, and so the graded family $\I$ is Noetherian by \cite[Theorem 2.1]{HHT07}.
\end{itemize}
Furthermore, it is easy to see that $\C(\I)$ is described by $x+y \ge 1, x \ge 0, y \ge 0$. Hence, it follows from Theorem \ref{thm.Mult}, that
$$e(\I)= 2!\text{co-vol}_2\left(\C(\I)\right)=1.$$
\end{example}

\begin{remark}
    \label{rmk.IdealsWithMonIntClosure}
    In a polynomial ring $R = \kk[x_1, \dots, x_d]$, NND ideals can be defined in similar way as in Definition \ref{def.NND}. However, in this setting, NND ideals are not characterized by having monomial integral closures. We call an ideal $I \subseteq R$, for which $\bar{I\ }$ is a monomial ideal, a \emph{weakly} NND ideal. It is desirable to study properties of weakly NND ideals.
\end{remark}


\end{document}